\documentclass[a4paper,12pt]{article}
\usepackage{amsfonts,amssymb,amsmath,amsthm}
\usepackage[dvips]{graphics}
\usepackage[dvips]{graphicx}
\usepackage[dvips]{color}

\def\num{\hspace{-2mm}{\bf }\hspace{2mm}}

\newtheorem{st}{Statement}[section]

\newtheorem{propo}[st]{Proposition}

\newtheorem{cor}[st]{Corollary}

\newtheorem{thm}[st]{Theorem}
\newtheorem{lemm}[st]{Lemma}

\def\max{{\rm max\,}}

\def\cP{{\cal P}}

\def\Fin{{\textbf S}}
\def\card{{\rm card\,}}

\def\Proof:{ \vspace{-1.5mm} {\noindent\it Proof.}}
\def\gm{\vspace{0mm}}
\def\gd{\vspace{0mm}}

\def\Box{\rule{1.5mm}{1.5mm}}

\begin{document}

\renewcommand{\thefootnote}{\fnsymbol{footnote}}

 \title{\gd \gd How high can Baumgartner's $\cal I$-ultrafilters lie in the P-hierarchy ?}
 \author{ Micha\l \ Machura and Andrzej Starosolski}
 \date{\today}
\maketitle

\gd
\begin{abstract}
Under CH we prove that for any tall P-ideal $\cal I$ on $\omega$ and for any ordinal $\gamma \leq \omega_1$ there is an ${\cal I}$-ultrafilter (in the sense of Baumgartner), which belongs to the class ${\cal P}_{\gamma}$ of P-hierarchy of ultrafilters. Since the class of ${\cal P}_2$ ultrafilters coincides with a class of P-points, out result generalize theorem of Fla\v{s}kov\'a, which states that there are ${\cal I}$-ultrafilters which are not P-points.  
\end{abstract}


\gd
\section{Introduction}

Baumgartner in the article \textit{Ultrafilters on $\omega$} (\cite{Baum}) introduced a notion of ${\cal I}$-ultrafilters: 
\vspace{5mm}

Let ${\cal I}$ be an ideal on $\omega$.
A filter on $\omega$ is an ${\cal I}$-ultrafilters, if and only if, for every function $f\in \omega^{\omega}$ there is a set $U\in u $ such that $f[U] \in {\cal I}$.
\vspace{5mm}

 This kind of ultrafilters was studied by large group of mathematician. We shall mention only the most important papers in this subject from our point of view: J. Brendle \cite{Brendle}, C. Laflamme \cite{Laf}, Shelah \cite{Shelah2} ,\cite{Shelah}, B\l aszczyk \cite{Blaszczyk}. The theory of $\cal I$-ultrafilters was developted by Fla\v{s}kov\'a in a series of articles and in her Ph.D thesis \cite{FlasDok}.  
\vspace{5mm}

In \cite{FlasDok} Fla\v{s}kov\'a  proved under CH that for every tall P-ideal ${\cal I}$ that contains all singletons, there is an ${\cal I}$-ultrafilters, which is not a $P$-point. Later she succeeded to replace CH by the assumption $\frak{p}=\frak{c}$ \cite{Flas1}.
\vspace{5mm}

Ultrafilters on $\omega$ may be classified with respect to
sequential contours of different ranks, that is, iterations of the
Fr\'{e}chet filter by contour operations. This way an
$\omega_1$-sequence $\{\cP_\alpha \}_{1\leq\alpha\leq\omega_1}$ of
pairwise disjoint classes of ultrafilters - the P-hierarchy - is
obtained, where P-points correspond to the class $\cP_2$, allowing
us to look at the P-hierarchy as the extension of notion of P-point. 
The following theorem was proved by Starosolski, see \cite{Star-P-hier} Proposition 2.1: 

\begin{propo}
An ultrafilter $u$ is a P-point if and only if $u$ belongs to the class ${\cal P}_2$ in P-hierarchy. 
\end{propo}

Many inmportant information about P-hierarchy may be found in \cite{Star-P-hier}.
For additional information regarding sequential cascades and
contours one can look at \cite{DolMyn}, \cite{DolStaWat},
\cite{Dol-multi}, \cite{Star-ff}. However the 
most important definitions and conventions shall be repeated below.

Since $P$-point correspond to ${\cal P}_2$ ultrafilter in P-hierarchy of ultrafilters (more about P-hierarchy one can find below), it would interesting to know to which classes of P-hierarchy can belong ${\cal I}$-ultrafilters. In this paper we shall show that it can be any class ${\cal P}_{\alpha}$. Let us introduce all necessary definitions and tools. 
\vspace{5mm}

The set of natural numbers (finite ordinal numbers) we denote $\omega$. The filter considered in this paper will be defined on infinite countable set (except one indicated case). This will be usually a set $\max V$ of maximal elements of a cascade $V$ (see definition of cascade below) and we will often identify it with $\omega$ without indication.  The following convention we be applied without mentioning it: 
\vspace{5mm}

\textit{Conventions:} If $u$ is a filter on $A \subset B$, then we identify $u$ with the
filter on $B$ for which $u$ is a filter-base. In particular we identify principal ultrafilter on $\{ v \}$ with principal ultrafilter generated on $\omega$ by $v$. If $\cal F$ is a filter base, then by $\langle \cal F \rangle$ we denote a filter generated by $\cal F$.
\vspace{5mm}

The {\it cascade} is a tree $V$ without infinite branches and with a
least element $\emptyset _V$. A cascade is $\it sequential$ if for
each non-maximal element of $V$ ($v \in V \setminus \max V$) the set
$v^{+V}$ of immediate successors of $v$ (in $V$) is countably
infinite. We write $v^+$ instead of $v^{+W}$ if it is known in which
cascade the successors of $v$ are considered. If $v \in V \setminus
\max V$, then the set $v^+$ (if infinite) may be endowed with an
order of the type $\omega$, and then by $(v_n)_{n \in \omega}$ we
denote the sequence of elements of $v^+$, and by $v_{nW}$ - the
$n$-th element of $v^{+W}$.
\vspace{5mm}

The {\it rank} of $v \in V$ ($r_V(v)$ or $r(v)$) is defined
inductively as follows: $r(v)=0$ if $v \in \max V$, and otherwise
$r(v)$ is the least ordinal greater than the ranks of all immediate
successors of $v$. The rank $r(V)$ of the cascade $V$ is, by
definition, the rank of $\emptyset_V$. If it is possible to order
all sets $v^+$ (for $v \in V \setminus \max V$)  so that for each $v
\in V \setminus \max V$ the sequence $(r(v_n)_{n<\omega})$ is
non-decreasing, then the cascade $V$ is {\it monotone}, and we fix
such an order on $V$ without indication.
\vspace{5mm}

For $v \in V$ we denote by $v^\uparrow$ a subcascade of
$V$ built by $v$ and all successors of $v$. 
We write $v^{\uparrow }$ 
instead of $v^{\uparrow V}$ 
if we know in which cascade the subcascade is included.
\vspace{5mm}

One may assume that cascade $V$ is a family of subset of infinite countable set ($\omega$) and the order on $V$ is inclusion. Indeed cascade $V$ is isomorphic to a cascade $\bar{V}$ such that: 
\begin{itemize}
\item  $\emptyset_{\bar{V}} = \omega$; 
\item  $\bar{v}^{+}$ is a partition of $\bar{v}$ for every $\bar{v}\in \bar{V}$: \\ $\bar{v} = \bigcup \{ \bar{w}: \bar{w} \in \bar{v}^{+} \} $ and elements of $\bar{v}^{+}$ are disjoint. 
\item  $\bar{v}$ is singleton for every $\bar{v} \in \max \bar{V}$.
\end{itemize}
An isomorphis $\bar{} : V\to \bar{V}$ is given by formula $\bar{v} = \max v^{\uparrow}$.
\vspace{0.5cm}



If $\mathbb{F}= \{{\cal F}_s: s \in S \}$ is a family of filters on
$X$ and if $\cal G$ is a filter on $S$, then the {\it contour of $\{
{\cal F}_s \}$ along $\cal G$} is defined by
$$\int_{\cal G} \mathbb{F} = \int_{{\cal G}} \{ {\cal F}_s  : s \in S \}  =
\bigcup_{G \in {\cal G}} \bigcap_{s \in G} {\cal F}_s.$$
\vspace{5mm}

Such a construction has been used by many authors (\cite{Fro},
\cite{Gri1}, \cite{Gri2}) and is also known as a sum (or as a limit)
of filters. 

Operation of sum of filters we apply to define \textit{contour of cascade}:
Fix a  cascade $V$. 
Let ${\cal G}(v)$ be a filter on $v^{+}$ for every $v \in V \setminus \max V$.
For $v \in \max V$ let ${\cal G}(v)$ be a trivial ultrafilter on a singleton $\{ v \}$ (we can treat it as principal ultrafilter on $\max v$ according to convention we assumed).
This way we have defined a function $v \longmapsto  {\cal G}(v)$. 
 We define contour of every sub-cascade $v^{\uparrow}   $  inductively with respect to rank of $v$:
$$\int^{\cal G} v^{\uparrow} = \{ \{ v \} \}$$ 
 for $v\in \max V$  (i.e. $\int^{\cal G} v^{\uparrow}$ is just  a trivial ultrafilter on singleton $\{ v \} $) ; 
$$\int^{\cal G} v^{\uparrow} =\int_{{\cal G}(v)} \left\{ \int^{\cal G} w^{\uparrow} : w\in v^{+} \right\} $$ for $v \in V \setminus \max v$.
Eventually we put 
$$ \int^{\cal G} V = \int^{\cal G} \emptyset_{V}.$$
Usually we shall assume that all the filters ${\cal G}(v)$ are Frechet (for $v\in V\setminus \max V$). In that case we shall write $\int V$ instead of $\int^{\cal G} V$.

Filters defined similar way were considered in \cite{Kat1}, \cite{Kat2}, \cite{Dagu}, also.

Let $V$ be a monotone sequential cascade and let $u=\int V$.
Then a {\it rank $r(u)$} of $u$ is, by definition, the rank of $V$.

It
was shown in \cite{DolStaWat} that if $\int V= \int W$, then $r(V) =
r(W)$.
\vspace{5mm}

We shall say that a set $F$ \textit{meshes} a contour ${\cal V}$ ($F \# {\cal V}$) if and only if ${\cal V} \cup \{ F\}$ has finite intersection property i.e can be extended to a filter. If $\omega \setminus F \in {\cal V}$, then we say that $F$ is \textit{residual} with respect to ${\cal V}$ .
\vspace{5mm}

Let us define ${\cP}_\alpha$
for $1 \leq {\alpha <\omega_1}$ on $\beta\omega$ (see
\cite{Star-P-hier}) as follows: $u \in {\cP}_\alpha$ if there is no
monotone sequential contour $C_{\alpha}$ of rank $\alpha$ such that
$C_{\alpha}  \subset u$, and for each $\beta$ in the range
$1\leq\beta < \alpha$ there exists a monotone sequential contour
$C_\beta$ of rank $\beta$ such that $C_\beta \subset u$. Moreover,
if for each $\alpha < \omega_1$ there exists a monotone sequential
contour $C_\alpha$ of rank $\alpha$ such that $C_\alpha \subset u$,
then we write $u \in {\cP}_{\omega_1}$.
\vspace{5mm}




Let us consider a monotone cascade $V$ and a monotone sequential
cascade $W$. We will say that $W$ is a sequential extension of $V$
if:

1) $V$ is a subcascade of cascade $W$,

2) if $v^{+V}$ is infinite, then $v^{+V} = v^{+W}$,

3) $r_V(v) = r_W(v)$ for each $v \in V$.

Obviously, a monotone cascade may  have many sequential extensions.

Notice that if $W$ is a sequential extension of $V$ and $U \subset
\max V $, then $U$ is residual for $V$ if and only if $U$ is
residual for $W$.
\vspace{5mm}

It cannot be proven in ZFC that all the classes ${\cal P}_{\alpha}$ are nonempty.
The following theorem was proved in \cite{Star-P-hier} Theorem 2.8:
\begin{thm}\num The following statements are equivalent:
\begin{enumerate}
\item P-points exist,
\item $\cP_\alpha$ classes are non-empty for each countable successor
$\alpha$,
\item There exists a countable  successor $\alpha>1$ such that the
class $\cP_\alpha$ is non-empty.
\end{enumerate}
\end{thm}

Starosolski has proved in \cite{Star-Top-P-Hier} Theorem 6.7  that: 
\begin{thm}\num Assuming CH every class ${\cal P}_{\alpha}$ is nonempty 
\end{thm}

The main theorem presented in this paper is on the one side an extension of Starosolski's result, but on the side based on it. 
\vspace{5mm}

Let us consider another technical notion which one could called a "restriction of a cascade".
Let $V$ be a monotone sequential cascade and let a set $H$ meshes the contour $\int V$. By $V^{\downarrow H}$ we denote a biggest monotone sequential cascade such that $V^{\downarrow H} \subset V$ and $\max V^{\downarrow H} \subset H$. 
It is easy to see that $H \in \int V^{\downarrow H}$. 
\vspace{5mm}

\section{Lemmas}

The following lemmas will be used in the prove of a main theorem. 

The first lemma is one of  lemmas proved in 
\cite{Star-Ord-V-P}  (see: Lemma 6.3 ):

 \begin{lemm}\num
Let $\alpha < \omega_1$ be a limit ordinal and let $({\cal V}_n : n<\omega)$ be
a sequence of monotone
sequential contours such that $r({\cal V}_n)<r({\cal
V}_{n+1})<\alpha$ for every $n$ and that $\bigcup_{n < \omega}{\cal V}_n$ has finite intersection property. 
Then there is no monotone sequential contour
$\cal{W}$ of rank $\alpha$ such that ${\cal{W}} \subset
\langle \bigcup_{n<\omega} {\cal V}_n \rangle$.
\end{lemm}

Since the paper with a prove of the above lemma is not published yet, the authors decided to included a prove at the end of this paper in a appendix. 

As a corollary we get: 

\begin{lemm}
Let $\alpha < \omega_1$ be a limit ordinal, let $({\cal V}_n)_{n<\omega}$ be an increasing ("$\subset$")
sequence of monotone sequential
contours, such that $r({\cal V}_n)<\alpha$ and let ${\cal F}$ be a countable family of
sets such that $\bigcup_{n<\omega} {\cal V}_n  \cup {\cal F}$ has finite intersection property. Then
$\langle \bigcup_{n<\omega}{\cal V}_n  \cup {\cal F}) \rangle$ do not contain any monotone sequential contour of rank $\alpha$.
\end{lemm}

\textit{Proof:} 
Assume that  $\cal F$ is finite. Let ${\cal W}_n = \{ U \cap \bigcap {\cal F} : U\in {\cal V}_n  $. It is easy to see that ${\cal W}_n$ is monotone sequential contour of the same rank as ${\cal V}_n$. Consider a sequence  $({\cal W}_n)$. By Lemma 2.1 the union $({\cal W}_n)$ do not contains contour of rank $\alpha$. 
 
Assume that $\cal F$ is infinite. Order $\cal F$ in $\omega$ type,
obtaining a sequence $(F_n)_{n<\omega}$. Next put $${\cal W}_n = \{ U \cap \bigcap_{i\leq n} F_i : U\in {\cal V}_n \}. $$ 
Consider a sequence  $({\cal W}_n : n<\omega )$ and use again Lemma 2.1 to show that the union $({\cal W}_n : n<\omega)$ do not contains contour of rank $\alpha$.  $\hfill$ $\Box$
\vspace{5mm}

The following lemma is a straightforward extension of the claim contained in the proof of \cite{Flas1} 
Theorem 3.2. and since a proof is almost identical to the quoted one,
we left it to the reader.

\begin{lemm}
Let $\cal I$ be a tall P-ideal that contains all singletons, let $\{ U_n : n<\omega \}$ be a pairwise disjoint sequence
of subsets of $\omega$, let $\{u_n : n<\omega \}$ be a sequence of $\cal I$-ultrafilters such that $U_n \in u_n$, finally let
$v$ be another one $\cal I$-ultrafilter. Then $\int_{v} \{ u_n : n<\omega \} $ is a $\cal I$-ultrafilter.
\end{lemm}

As immediate consequence we get 

\begin{lemm}
If $V$ is monotone sequential cascade, 
${\cal G}(v)$ is an P-point and $\cal I$-ultrafilter for each $v \in V \setminus \max V$
and ${\cal G}(v)$ is a trivial ultrafilter on a singleton $\{ v \}$ for $v\in \max V$, then 
$\int^{\cal G} V$ is an $\cal I$-ultrafilter.
\end{lemm}

Similar lemma as above one can formulate for ultrafilters in certain class in P-hierarchy instead of $\cal I$-ultrafilters, see \cite{Star-P-hier} Theorem 2.5:

\begin{thm} Let $\gamma$ be an ordinal. 
Let $V$ be a monotone sequential cascade of rank $\gamma$, let $G(v)$ be a principal ultrafilter on $\{ v\} $ for $v \in \max V$, and let $G(v)$ be a P-point on $v^+$ for
$v \in V\setminus \max V$. Then $\int^GV \in P_{\gamma+1}$.
\end{thm}

\begin{cor}
If $V$ is monotone sequential cascade, 
${\cal G}(v)$ is an ultrafilter from the class $\cal P_{\gamma}$ for each $v \in V \setminus \max V$
and ${\cal G}(v)$ is a trivial ultrafilter on a singleton $\{ v \}$ for $v\in \max V$, then 
$\int^{\cal G} V$ belongs to the class $\cal P_{\gamma}$.
\end{cor}

\section{Main result}

In this section we shall present main result of the paper. 

\begin{thm} (CH) Let $\cal I$ be a tall P-ideal that contain all singletons, 
and let $\gamma \leq \omega_1$ be an ordinal.
Then there exists an $\cal I$-ultrafilter $u$ which belongs to $\cP_{\gamma}$.
\end{thm}
\Proof:
We shall split proof into five cases: $\gamma=1$, $\gamma=2$, $\gamma >2$ is a succesor ordinal (the main step), $\gamma <\omega_1$ is limit ordinal, $\gamma= \omega_1$. 
\vspace{0.5cm}

\noindent \textbf{ Step 0:} $\gamma =1$ is clear, image of singleton
($\cP_1$ is a class of principal ultrafilters) is a singleton, so belongs to $\cal I$.
\vspace{5mm}

\noindent \textbf{ Step 1:}  for $\gamma = 2$. 

We order all contours of rank 2 and all functions $\omega \rightarrow \omega$ in $\omega_1$-sequences
$({\cal W}_\alpha)_{\alpha<\omega_1}$, $(f_\alpha)_{\alpha<\omega_1}$ respectively . 
By transfinite induction, for $\alpha < \omega_1$ we build countable generated filters ${\cal F}_\alpha$ together with their decreasing basis $( F^n_{\alpha})_{n<\omega}$,
such that:
\begin{enumerate} 

\item ${\cal F}_0$ is a Frechet filter;

\item for each $\alpha < \omega_1$, the sequence $({ F}_\alpha^n)_{n<\omega}$is 
strictly decreasing base of ${\cal F}_{\alpha}$; 

\item ${\cal F}_\alpha \subset {\cal F}_\beta$ for $\alpha<\beta$;

\item ${\cal F}_\alpha = \bigcup_{\beta<\alpha}{\cal F}_\beta$ for $\alpha$ limit ordinal;

\item for each $\alpha<\omega_1$ there is $ F \in {\cal F}_{\alpha+1}$ such that $f_\alpha[F]\in {\cal I}$;

\item for each $\alpha<\omega_1$ there is $F\in{\cal F}_{\alpha+1}$ such that a complement of $F$ belongs to ${\cal W}_\alpha$.
\end{enumerate}

Suppose that ${\cal F}_\alpha$ is already define, we will show how to build ${\cal F}_{\alpha+1}$.
Since $F_\alpha^n$ is strictly decreasing one can pick $x_n \in F^n_{\alpha} \setminus F^{n+1}_{\alpha}$ for every $n<\omega$. Put $T= \{ x_n : n<\omega \}$. The are two possibilities:  

If $f_\alpha[T]$
is finite then there is $j \in f_\alpha[T]$ such that
a preimage $f_\alpha^{-1}[{j}]$ intersect infinite many  of $F_\alpha^n \setminus F_\alpha^{n+1}$. In this case put $G=f_\alpha^{-1}[{j}]$.

If $f_\alpha[T]$ is infinite, then since $\cal I$ is tall there is $I \in {\cal I}$ such that $I \subset f_\alpha[T]$.
This time put $G=f_\alpha^{-1}[I]$.
\vspace{0.5cm}

Notice that $\{F_\alpha^n : n<\omega \}\cup \{G_{\alpha} \}$ has finite intersection property and is countable. A subbase of any sequential contour of rank 2 has cardinality at least ${\frak d} >\aleph_0$, thus none of them one is contained in $\{F_\alpha^n : n<\omega \}\cup \{G_{\alpha} \}$.
This means that there is a set $A_{\alpha}$ such that its complement belongs to  ${\cal W}_{\alpha}$ and a family  $\{F_\alpha^n : n<\omega \}\cup \{G_{\alpha}, A_{\alpha} \}$
has finite intersection property. 
Order $ \{F_\alpha^n : n<\omega \}
\cup\{G_\alpha\} \cup \{A_\alpha\}$ in $\omega$ type, obtaining a sequence $(\tilde{F}_{\alpha+1}^n : n<\omega )$. 
Put $F_{\alpha +1}^n = \bigcap_{m\leq n}\tilde{F}_{\alpha +1}^n$ to get decreasing sequence and let ${\cal F}_{\alpha + 1}=
\langle \{F_{\alpha+1}^n: n < \omega\} \rangle$.
\vspace{0.5cm}

Take any ultrafilter $u$ that extends
$\bigcup_{\alpha<\omega_1}{\cal F}_\alpha$.
By condition 5) $u$ is an $\cal I$-ultrafilter, by condition 6)
$u$ do not contain any monotone sequential contour of rank $2$. Since by condition 1) $u$
 contains a Frechet filter it is not principal. Thus $u$ is a P-point.
(Note that on this step we do not use an assumption, that $\cal I$ is a P-ideal.)

\hspace{5mm}

\textbf{Step 2:} $\gamma$ is an arbitrary successor ordinal such that $2<\gamma<\omega_1$.
Let $V$ ba an arbitrary monotone sequential cascade of rank $\gamma - 1$. Let $V \in v \longmapsto {\cal G}(v)$ be a function such that: 

1) ${\cal G}(v)$ is an P-point and $\cal I$-ultrafilter for each $v \in V \setminus \max V$
(such ultrafilters exists by step 1)

2) ${\cal G}(v)$ be a trivial ultrafilter on a singleton $\{ v \}$ for $v\in \max V$.

 Lemma 2.5 guarantee that 
$\int^{\cal G} V \in \cP_\gamma$. 
whilst Lemma 2.4 ensures us that  $\int^{\cal G} V$ is an $\cal I$-ultrafilter.

So we are done for successor $\gamma$.

\hspace{5mm}

\textbf{Step 3:} for limit $\gamma <\omega_1$.
The proof in this case is base on the same idea as step 1, but it is more sophisticated and technical. 

Let $({\cal V}_n)_{n<\omega}$ be an increasing ("$\subset$")
sequence of monotone sequential
contours, such that their ranks $r({\cal V}_n)$ are smaller than $\gamma$ but converging to $\gamma$. For each $n< \omega$ denote by $V_n$ a
(fixed) monotone sequential cascade
such that $\int V_n = {\cal V}_n$.
Let $\{{\cal W}_{\alpha}, \alpha<\omega_1\}$ be an enumeration of all monotone sequential contours
of rank $\gamma$.
Let $\omega^\omega = \{f_\alpha: \alpha< \omega_1\}$.
\hspace{0.5cm}

By transfinite induction, for $\alpha < \omega_1$ we build filters ${\cal F}_\alpha$ together with their decreasing basis $({ F}_\alpha^n)_{n<\omega}$,
such that:
\begin{enumerate}
\item ${\cal F}_0$ is a Frechet filter;

\item for each $\alpha < \omega_1$ $({ F}_\alpha^n)_{n<\omega}$ is a
strictly decreasing base of ${\cal F}_\alpha$;

\item ${\cal F}_\alpha \subset {\cal F}_\beta$ for $\alpha<\beta$;

\item ${\cal F}_\alpha = \bigcup_{\beta<\alpha}{\cal F}_\beta$ for $\alpha$ limit ordinal;

\item $\bigcup_{i<\omega}{\cal V}_i \cup \bigcup_{\alpha<\omega_1}{\cal F}_\alpha$
has finite intersection property;

\item for each $\alpha<\omega_1$ there is $ F \in {\cal F}_{\alpha+1}$ such that $f_\alpha[F]\in {\cal I}$;

\item for each $\alpha<\omega_1$ there is $F\in{\cal F}_{\alpha+1}$ such that the complement of $F$ belongs to ${\cal W}_{\alpha}$.
\end{enumerate}

Suppose that ${\cal F}_\alpha$ is already define, we will show how to build ${\cal F}_{\alpha+1}$.
This shall be done in five substeps. 
First for each ${\cal V}_n$ and each $F^i_{\alpha}$ we shall find $H_{n,i}$ such that ${\cal V}_{n} \cup \{ F^i_{\alpha} , H_{n,i} \}$ has finite intersection property and $f_{\alpha}[H_{n,i}] \in {\cal I}$. 
Next we shall replace all the sets $H_{n,i}$ by one set  $H_n$ such that 
${\cal V}_n \cup {\cal  F}_{\alpha} \cup \{ H_{n} \}$ has finite intersection property and $f_{\alpha}[H_{n}] \in {\cal I}$.
On the third step one has to replace all the sets $H_{n}$ by one set $G_{\alpha}$ such that 
$\bigcup_{m<\omega} {\cal V}_n \cup {\cal  F}_{\alpha} \cup \{ G_{\alpha} \}$ has finite intersection property and $f_{\alpha}[G_{\alpha}] \in {\cal I}$.
The set $G_{\alpha}$  take care on all the contours ${\cal V}_n$. Adding it as generator  to $F_{\alpha +1}$ will ensure preservation of conditions 5 and 6. 
On the fourth step will take care on condition 7 by adding set $A_{\alpha}$ to the list of generators of $F_{\alpha +1}$.
The last thing is to define decreasing base of a filter ${\cal F}_{\alpha +1}$ and a filter itself. 

\vspace{0.5cm}

\noindent \textit{Substep i)}  Fix $n$ and $i$. Let us introduce an axillary definition.
\vspace{0.5cm}

\noindent
\textit{Definition:} Fix a monotone sequential cascade $V$, a set $F$ and  a function $f \in \omega^{\omega}$. For each $v\in V$, we write $U\in \Fin(v)$ if
\begin{enumerate}
\item  $U \subset \max v^\uparrow$;

\item $(U \cap F) \# \int v^\uparrow$;

\item $\card(f[U \cap F])=1$.
\end{enumerate}

We following claim is crucial:

\begin{propo}
One that one of the following possibilities holds:
\vspace{5mm} 
 
A) $\Fin(\emptyset_{V})\not=\emptyset$;
\vspace{5mm}

B) there is an antichain (with respect to the order of a cascade) $\mathbb{A}\subset V$ such that:
\begin{enumerate}
\item  $\Fin(v) = \emptyset$ for all $v \in \mathbb{A}$,

\item $ \left( \bigcup \{  \max w^\uparrow   : w\in v^+, \Fin(w)\not= \emptyset  \} \right)  \# \int v^\uparrow$ for all $v \in \mathbb{A}$,

\item $\left( \bigcup \{
\max v^\uparrow :   v \in {\mathbb{A}}  \} \right) \# \int V$.
\end{enumerate}
\end{propo}

\textit{Proof of the proposition.} 
First notice, that in definition of $\Fin$ one can replace cardinality one by finite in condition 3), and that the replacement do not influence non-emptyness of $\Fin(v)$.

The proof is inductive by the rank of cascade $V$.
\vspace{5mm}

\noindent \textit{First step:} $r(V)=1$.  If case A holds, then we are done, so without loss of generality
$f (U \cap F)$ is infinite for each $U \cap F \in \max V$ such that $U \# \int V$.
But since $r(V)=1$, thus $\card(f(\max w \cap F)) \leq 1$, for each $w\in v^+$. And since $F \# \int V$ thus 
$$\left( \bigcup \left\{  (\max w \cap F) : w \in v^+, \card(f(\max w \cap F))=1   \right\} \right)   \# \int V   , $$
We put $\mathbb{A}= \{ \emptyset_V \} $ and see that case B holds.
\vspace{5mm}

\noindent \textit{Inductive step :} 
Suppose that the proposition is true for each $\beta<\alpha<\omega_1$. So take $V$ that $r(V)=\alpha$.
Again if case A holds, then we are done, so without loss of generality assume that
$f (U \cap F)$ is infinite for each $U \cap F \subset \max V$ such that $U \# \int V$.
By inductive assumption, for each successor $w$ of $\emptyset_{V}$ either case A holds for cascade $w^{\uparrow}$ 
either case B holds for for cascade $w^{\uparrow}$ . 

Split the set $\emptyset_{V}^{+}$ of immediate succesors of $\emptyset_{V}$ into to subsets: 
 $$V^A= \left\{ w \in \emptyset_V: \mbox{ case A holds } \right\}, \ \
V^B= \left\{ w \in \emptyset_V:  \mbox{case B holds }  \right\}.$$ 
Since $F \# \int V$, we have two possibilities: 

\noindent
$$\left( \bigcup_{w\in V^A}(\max w^\uparrow \cap F) \right) \ \# \int V \mbox{ or }  
 \left( \bigcup_{w\in V^B}(\max w^\uparrow \cap F) \right) \ \# \int V .$$ 

In the first case $\mathbb{A}= \{  \emptyset_V \} $ we was looking for.

In the second case, for each $w\in V^B$ there is a claimed (by inductive assumptions)
antichain ${\mathbb{A}}_{w}$ in $w^\uparrow$. Put ${\mathbb{A}} = \bigcup_{w \in V^B}{\mathbb{A}}_{w}$. This finishes proof of the proposition. \hfill $\Box$ 
\vspace{5mm}

We can come back to the main proof.

We aplly a proposition to  cascade $V_n$, set $F^i_{\alpha}$ and a function $f_{\alpha}$
In the case A we take any $U \in \Fin(\emptyset_{V_n})$ and  denote it by $H_{n,i}$.
\vspace{5mm}

In the case B for any $v \in \mathbb{A}$ we fix $U_{w}\in \Fin(w)$ for every $w\in v^+$ for which $\Fin(w)\not= \emptyset$;
for all the other $w\in V_n$ let $U_{w}=\emptyset$.
For $v \in \mathbb{A}$ consider $T_{v}=\bigcup_{w \in v^+} U_{w}$, and 
notice that $f_\alpha[T_v]$ is infinite. 
Since $\cal I$ is tall there is an infinite $I_v\in {\cal I}$ such that $I_v \subset  f_\alpha[T_v]$. 
Since $\cal I$ is an P-ideal, there is infinite $I_{n,i}\in {\cal I}$ such that $I_v \setminus I_{n,i}$
is finite for all $v \in \mathbb{A}$. Put $H_{n,i}=f^{-1}[I_{n,i}]$.
\vspace{0.5cm}

\noindent \textit{Substep ii)}
Now we will show how to replace sets $H_{n,i}$  by one set $H_n$. Consider two possibilities:

C) there is an infinite $K \subset \omega$ that $f_\alpha[H_{n,i}]$ is infinite for each $i\in K$;

D) there is an infinite $K \subset \omega$ that $f_\alpha[H_{n,i}]$ is a singleton for each $i\in K$.

\noindent In both cases since since $(F_\alpha^i)_{i<\omega}$ is decreasing,
without loss of generality we may assume that $K = \omega$.
\vspace{5mm}

In the case C, since $\cal I$ is an P-ideal,
there is infinite $I_{n}\in {\cal I}$ such that $I_n \setminus I_{n,i}$ is finite for each $i<\omega$.
Put $H_n = f_\alpha^{-1}[I_n]$.
\vspace{5mm}

In the case D we have two sub-cases: 

If $f_\alpha[\bigcup_{i<\omega}H_{n,i}]$ is infinite, then since $\cal I$ is tall, there is an infinite $I_n\in {\cal I}$ such that
$I_n \subset f_\alpha[\bigcup_{i<\omega}H_{n,i}$, and we put $H_n = f_\alpha^{-1}[I_n]$.

If $f_\alpha[\bigcup_{i<\omega}H_{n,i}]$ is finite; then there is $j \in f_\alpha[\bigcup_{i<\omega}H_{n,i}]$ that $f^{-1}_\alpha[\{j\}] =H_{n,i}$
for infinite many $i$'s, and we put $H_n = f_\alpha^{-1}[\{j\}]$.
\vspace{5mm}

Clearly, in both  cases  ${\cal V}_n \cup {\cal F}_\alpha  \cup \{H_n\}$ has finite intersection property
and $f_\alpha[H_n] \in {\cal I}$. 

\vspace{5mm}

\noindent \textit{Substep iii)}
On this step we have to find set $G_{\alpha}$ which can replace each $H_n$.
 We have shown that for each $n$ there is a set $H_n$  such that
and $f_{\alpha}[H_n] \in {\cal I}$.  In fact we got a little bit more: either $f_{\alpha}[H_n]$ is infinite but belongs to ${\cal I}$, either $f_{\alpha}[H_n]$ is a singleton. We set 
$$ S= \{ n<\omega : ( \exists R_n ) \ : {\cal V}_n \cup {\cal F}_\alpha  \cup \{ R_n\} \mbox{ has f.i.p. and } f_{\alpha}[R_n] \mbox{ is singleton } \}$$
It could happen that  $f_{\alpha}[H_n]$ is infinite but $n\in S$ and for some $R_n$ as above an image $f_{\alpha}[R_n]$ is singleton.
It this case we replace $H_n$ by any $R_n$. 
For $n \in \omega \setminus S $ we leave $H_n$ unchanged.
Once again proof splits into two cases: either $S$ is infinite, either it is finite. 
\vspace{5mm}

For infinite  $S$: Without loss of generality (since $({\cal V}_n)$ is increasing) we may assume that $S=\omega$ i.e. 
$f_\alpha[H_n]$ is a singleton for each $n<\omega$.
\vspace{5mm}

If $f_\alpha[\bigcup_{n<\omega} H_n ]$ is finite, then there is $j \in f_\alpha[\bigcup_{n<\omega} H_n ]$ such that 
$f_\alpha[H_n]=\{ j \}$ for infinite many $n$. Since ${\cal V}_n$ is increasing and $(F_\alpha^n)$ is decreasing,
a family $ \bigcup_{n<\omega} {\cal V}_n \cup {\cal F}_\alpha \cup f_\alpha^{-1}[\{j\}])$
has finite intersection property. Put $G_\alpha=f_\alpha^{-1}[\{j\}]$.
\vspace{5mm}

If $f_\alpha[\bigcup_{n<\omega} H_n ]$ is infinite, then, since $\cal I$ is tall,
there is infinite $I_\alpha\in {\cal I}$ such that $I_\alpha\subset f_\alpha[\bigcup_{n<\omega} H_n ]$.
Since ${\cal V}_n$ is increasing and $(F_\alpha^n)$ is decreasing,
a family $\bigcup_{n<\omega} {\cal V}_n \cup {\cal F}_\alpha \cup f_\alpha^{-1}[I_\alpha]
$
has finite intersection property. Put $G_\alpha=f_\alpha^{-1}[I_\alpha]$.
\vspace{5mm}

For finite  $S$: Without loss of generality
(since $({\cal V}_n)$ is increasing) we may assume that $S=\emptyset$ i.e. $f_\alpha[H_n]$ is infinite
for each $n<\omega$.

Since $\cal I$ is a P-ideal, and $f_\alpha[H_n] \in {\cal I}$,  there is $I_\alpha \in {\cal I}$
such that $f_\alpha[H_n] \setminus I_{\alpha}$ is finite for each $n<\omega$.

Since the sequence $({\cal V}_n)$ is increasing, we have two possibilities: 
either $f^{-1}_\alpha[I_\alpha]\# {\cal V}_n$ for all $n<\omega$; either  $ \neg f^{-1}_\alpha[I_\alpha]\#{\cal V}_{n}$ for almost every $n$.
The second possibility cannot happen by the definition of $\omega\setminus S$.
Put $G_\alpha = f^{-1}_\alpha[I_\alpha]$. It is easy to see that 
a family 
$\bigcup_{n<\omega}{\cal V}_n \cup {\cal F}_\alpha \cup \{G_\alpha\}$
 has finite intersection property.
\vspace{5mm}

\noindent \textit{Substep iv)} 
 Since the family ${\cal F}_\alpha \cup \{G_\alpha\}$  is countable, thus by Lemma 2.2 
 there exists $A_\alpha$ residual for the contour ${\cal W}_{\alpha}$ and such that a family
 $\bigcup_{n<\omega}{\cal V}_n \cup {\cal F}_\alpha
\cup\{G_\alpha , A_\alpha\}$ has finite intersection property.
\vspace{5mm}

\noindent \textit{Substep v)}
Order $ {\cal F}_\alpha
\cup\{G_\alpha\} \cup \{A_\alpha\}$ in type $\omega$, obtaining a sequence \\ $(\tilde{F}_\alpha^n : n<\omega )$.
Put $F_\alpha^n = \bigcap_{m\leq n}\tilde{F}_{\alpha}^n$ to get decreasing sequence and let ${\cal F}_{\alpha + 1}=
\langle \{F_{\alpha+1}^n: n < \omega\} \rangle$.
\vspace{0.5cm}

Take any ultrafilter $u$ that extends
$\bigcup_{n<\omega}{\cal V}_n \cup \bigcup_{\alpha<\omega_1}{\cal F}_\alpha$.
By condition 5) $u$ is an $\cal I$-ultrafilter, by condition 6)
$u$ do not contain any monotone sequential contour of rank $\gamma$ which jointly with 
$\bigcup {\cal V}_n \subset u$ give us $u \in {\cal P}_\gamma$.

So the proof is done also for limit $\gamma$.  
\vspace{5mm}

\textbf{Step 4:} $\gamma = \omega_1$. We will show a little more i.e. that there is a supercontour which
is an $\cal I$-ultrafilter.

Again we list $^\omega\omega=\{f_\alpha :   \alpha<\omega_1 \} $,
and we also list all pair (set and its complement) in the $\omega_1$-sequence of pairs
$(A_\alpha, \omega \setminus A_\alpha)$  that way that each set appears in the sequence only ones: or a set $A_{\alpha}$ or as complement $\omega \setminus A_{\alpha}$.

We will build an $\omega_1$ sequence $(V_\alpha : \alpha<\omega_1 )$ of monotone sequential cascades 
such that

\begin{enumerate}
\item $\int V_\beta \subset \int V_\alpha$ for each $\beta<\alpha<\omega_1$.

\item $r(V_\alpha)=\alpha$ for every $\alpha < \omega_1$ ;

\item $\max v_{\alpha} = \omega$ for every $\alpha<\omega_1$; 

\item there exist $U \in \int V_{\alpha+1}$ such that $f_\alpha[U] \in {\cal I}$

\item $A_\alpha \in \int V_{\alpha+1}$ or $\omega \setminus A_\alpha \in \int V_{\alpha+1}$.
\end{enumerate}

Define $V_1$ as an arbitrary (fixed) monotone sequential cascade of rank 1. Suppose that we already defined
cascades $V_\beta$ for all $\beta<\alpha<\omega_1$.
\vspace{5mm}

Case 1) $\alpha = \beta+1$ is a successor. Take $V_{\beta}$, by step 3 there is a set $H_{\alpha}$ suth that  $H_\alpha \# \int V_{\beta}$
and $f_\delta[H_\delta]\in {\cal I}$.
Consider a cascade $V_{\beta}^{\downarrow  H_\alpha}$; this is a monotone sequential cascade of
rank $\beta$. By the proof of Theorem 4.6 from \cite{DolStaWat} 
there is a monotone sequential cascade $\tilde{V}_{\alpha} $ of rank $\alpha$ such that
$\int V_{\beta}^{\downarrow  H_\alpha} \subset \int \tilde{V}_{\alpha}$. At least one of the elements of a pair
$(A_\alpha, \omega \setminus A_\alpha)$ mashes $\int \tilde{V_\alpha}$, denote it by   $B_\alpha$. Now let $V_\alpha
= \tilde{V}_{\alpha}^{\downarrow B_\alpha}$.
\vspace{5mm}

Case 2) $\alpha$ is limit. Let $V_\alpha$ be any monotone sequential cascade of rank $\alpha$ such that
$\int V_\beta \subset \int V_\alpha$ for each $\beta< \alpha$. Such a
cascade was constructed in the proof of Theorem 4.6 in \cite{DolStaWat}.  
\vspace{5mm}

Now it suffice to take $u=\bigcup_{\alpha<\omega_1} \int V_\alpha$. 
By construction $u$ has a finite intersection property and is a supercontour,
by 4) $u$ is an ultrafilter and by 3) $u$ is an $\cal I$-ultrafilter. \hfill $\Box$
\vspace{5mm}

The assumption that an ideal ${\cal I}$ is tall is essential: Fla\v{s}kov\'a has proved in \cite{Flas1} Proposition 2.2,  that if ${\cal I}$ is not tall, then there is no $\cal I$-ultrafilters. One can easily see, that an  ideal ${\cal I}$  has to contains all singletons, also.

\section{Appendix}

In this appendix we shall prove Lemma 2.1. The main tools we  use is an operation of decreasing the rank of cascade described below. Let us introduce axillary notion: 
\vspace{5mm}

Let $V$ be a cascade and let $x_0, x_1, x_2, \ldots $ be immediate succesors of $\emptyset_V$. We denote a sub-cascades $x_i^{\uparrow }$ by $V^{(i)}$. Similarly for $V^{(i)}$ if $x_{i0}, x_{i1}, x_{i2}, \ldots $ are immediate succesors of $\emptyset_{V^{(i)}}=x_i$ then we denote sub-cascades $x_{ij}^{\uparrow }$ by $V^{(i)(j)}$.
\vspace{0.5cm}

We say that a cascade $V$ is built by destruction of nods of rank 1 in a cascade $W$ iff 
\begin{itemize}
\item all elements of rank 1 are removed from $W$: \\ $V= W\setminus \{ v\in W : r_W(v)=1 \}$;
\item immediate succesors of elements which had rank 2 are succesors of their former succesors: if $r_W(v)=2$ then
\end{itemize}
$$v^{+V} = \bigcup \left\{ w^{+W} : w\in v^{+W}  \right\}.  $$

Observe that if $r(W)$ is finite then $r(V) = r(W)-1$.

Assume that there we are given a cascade of rank $\alpha$ and an ordinal $\beta <\alpha$. We shall describe a operation of decreasing of rank of a cascade $W$. 
The construction is inductive: 
\vspace{0.5cm}

\textit{$\alpha$ is finite:} 

We can decrease rank of $W$ from $\alpha$ to $\beta$ by applying $\alpha - \beta$ times an operation of destructing nods of rank 1. 
\vspace{0.5cm}

\textit{$\alpha$ is infinite:}
\vspace{0.5cm} 

\textit{$\beta = \bar{\beta} +1$  is succesor ordinal} and we are able to decrease a rank of any cascade of rank smaller than $\alpha$. Let $r(W)=\alpha$. Consider cascades $W^{(i)}$ for $i<\omega$. Of course $r(W^{(i)}) <\alpha$ for every $i$ and one can decrease their ranks to $\bar{\beta}$. Let $V^{(i)}$ be cascades obtained from $W^{(i)}$ by decrasing rank: $r(V^{(i)})= \bar{\beta}$, and  let $V$ be a cascade obtained  by gluing cascases $V^{(i)}$ together. Thus $r(V)=\bar{\beta}+ 1=\beta$.
\vspace{0.5cm}

\textit{$\beta$  is limit ordinal} and we are able to decrease a rank of any cascade of rank smaller than $\alpha$. Let $r(W)=\alpha$. Consider cascades $W^{(i)}$ for $i<\omega$ and a sequence of ordinal $(\beta_i)_{i<\omega}$ increasing to $\beta$. Of course $r(W^{(i)}) <\alpha$ for every $i$ and one can decrease ranks of every $W^{(i)}$ to $\beta_i$. Let $V^{(i)}$ be cascades obtained from $W^{(i)}$ by decrasing rank: $r(V^{(i)})= \beta_i$, and  let $V$ be a cascade obtained  by gluing cascases $V^{(i)}$ together. Thus $r(V)=\beta$.
\vspace{0.5cm}

Observe that above desribed decreasing of rank of cascade $W$ does not change $\max W$.
If a cascade $W$ is obtained from $V$ by decreasing rank, then we write $W \triangleleft V$. 
Trivially $V \triangleleft V$ for every $V$. 
\vspace{0.5cm}

We shall make use of the following theorem (see: \cite{Dol-multi} ) :

\begin{thm}[Dolecki]  
If $({\cal V}_n)_{n<\omega}$ is a sequence of monotone sequential contours of rank less than $\alpha$ and 
$\bigcup_{n<\omega} {\cal V}_n$ has finite intersection property, then there is no monotone sequential countour ${\cal W}$ of rank $\alpha +1$  such that ${\cal W} \subset \langle \bigcup_{n<\omega} {\cal V}_n \rangle $.
\end{thm}

Before we prove Lemma 2.1 we shall prove a following technical claim

\begin{lemm}
Let $V$ be a cascade of rank $\alpha$, $W$ be cascade obtained from $V$ by decreasing rank of $V$ to $\beta <\alpha$ and let $\beta< \gamma < \alpha$. Then there is a cascade $T$ of rank $\gamma$ such that $W\triangleleft  T \triangleleft  V$.  
\end{lemm}

\textit{Proof:} 
The proof is inductive on triples $(\beta, \alpha, \gamma)$ where $\beta \leq \gamma \leq \alpha$ and  ordered lexographically.
Assume that for $(\beta', \alpha', \gamma') < (\beta, \alpha, \gamma)$ lemma has been proved. 
For $\gamma= \beta$ there is nothing to prove, so assume that $\beta <\gamma$.  
\vspace{5mm}

Observe that if $v\in W$ is an element of rank $1$ in $W$, then  its succesors are maximal elements of cascade $V$: 
$$ v^{+W} = \max v^{\uparrow V}. $$
Consider two cases.
\vspace{0.5cm}

\textit{ $\gamma$ is  a limit ordinal:} Denote by $x_1, x_2, \ldots $ succesors of $\emptyset_V $. Recall that  $ \emptyset_V = \emptyset_W$ and $ \emptyset^{+}_V = \emptyset^{+}_W$.
Fix increasing sequence $(\gamma_n : n<\omega)$ converging to $\gamma$ such that $\gamma_n > r_V(x_n)$ for every $n$.
We have 
$$ x_n^{\uparrow W}  \triangleleft x_n^{\uparrow V} \mbox{ and  }  r(x_n^{\uparrow W}) <\beta \  , \  r(x_n^{\uparrow V}) <alpha . $$
Using inductive hypothesis one can find $T^{(n)}$ such that $r(T^{(n)}) = \gamma_n $ and
$$ x_n^{\uparrow W} \triangleleft  T^{(n)}    \triangleleft x_n^{\uparrow V} . $$
Let $T$ be obtained by gluing $T^{(n)}$.

\vspace{5mm}
\textit{ $\gamma = \delta+1 $ is  a succesor ordinal:} Proceed similarly, by assume that a sequence  
$(\gamma_n : n<\omega) $ is contantly equal $\delta$. 
\hfill $\Box$ 
\vspace{5mm}

Now we can turn out attention to the proof of Lemma 2.1.
\vspace{0.5cm}

\textit{Proof of Lemma 2.1:}

Assume that there exists a contour ${\cal W}$ of rank $\alpha$ such that 
${\cal W} \subset \langle \bigcup_{n<\omega} {\cal V}_n \rangle$. 
We build a cascade $W$ and a sequence of cascades $(W_n)_{n<\omega}$ such that:
\begin{itemize}
\item $\int W ={\cal W}$;
\item $W_m \triangleleft W_{m+1}$ for every $m$;
\item $W_m \triangleleft W$ for every $m$;
\item $W_m$ is obtained by decreasing of rank of $W$ (with cutted several branches the way not influencing contour) to $\alpha_m +3$;
\item if $r(W_m^{(i)}) = \alpha_m+2$ then $r(W_m^{(i)(j)}) = \alpha_m+1$ for each $j$;
\item if $r(W_m^{(i)}) < \alpha_m+2$ then $W_m^{(i)} =  W_{m-1}^{(i)}$.
\end{itemize}
\vspace{0.5cm}

Fix any cascade $\bar{W}$ such that $\int \bar{W} = {\cal W}$.
Let $\bar{W}_m$ be a cascade obtained from $\bar{W}$ by cutting every branch $\bar{W}^{(i)}$ of rank smaller than $\alpha_m +2$ and every (sub-)branch  
$\bar{W}^{(i)(j)}$ of rank smaller than $\alpha_m +1$. Observe that we cut only finite many branches $\bar{W}^{(i)}$ and for the other $\bar{W}^{(i)}$ only finite many branches 
$\bar{W}^{(i)(j)}$. Thus $\int \bar{W}_m = \int \bar{W} = {\cal W}$ for every $m$.   
\vspace{0.5cm}

Let $W_1$ be a cascade obtained from $\bar{W}_1$ by decreasing ranks of $W_1^{(i)(j)} $  to $\alpha_1 +1$ and let $W=\bar{W}$. Thus 
$W_1 \triangleleft W $.  Assume that $W_1 \triangleleft W_2 \triangleleft \ldots \triangleleft  W_m$ have been defined such that $W_l \triangleleft W$ and $r(W_l^{(i)(j)})= \alpha_l +1$ (thus $r(W_l)= \alpha_l +3$)   for every $l\leq m$. We use Lemma 4.2 to cascades    $W_m^{(i)(j)}$ and $W^{(i)(j)}$ to define $W_{m+1}^{(i)(j)}$ of rank $\alpha_{m+1}+1$. Gluing $W_{m+1}^{(i)(j)}$ we obtain first $W_{m+1}^{(i)}$ and next gluing them $W_{m+1}$ such that $W_m \triangleleft W_{m+1} \triangleleft  W$ and $r(W_{m+1})= \alpha_{m+1} +3$.  
\vspace{0.5cm}

Next we build a decreasing sequence $(U_n)_{n<\omega}$ such that:
\begin{enumerate}
\item $U_n \in \int W_n$; 
\item $U_n \notin \langle \bigcup_{i\leq n} {\cal V}_i  \rangle $; 
\item $U_n \cap \max \bar{W}_{n+1} = U_{n+1} \cap \max \bar{W}_{n+1}$ for each $n$;
\item $ U_n  \cap \max W^{(k)} \in \int W^{(k)}$ for each $k$. 
\end{enumerate}
The last condition can follows difficulties in a construction. Therefore we shall replace in a construction cascades $W_n$ by modified cascades $\widetilde{W}_n$, such that it will be not necessary to take care on condition (4). Fix $n$. We define cascade $\widetilde{W}_n$ as follows:
We remove from original $W_n$ succesors of $\emptyset_{W_n}$ and we treat succesors of removed elements as new succesors of $\emptyset_{\widetilde{W}_n} = \emptyset_{W_n}$. Formally: 
$$ \emptyset_{\widetilde{W}_n}^{+} = \bigcup \left\{ w^{+} : w \in  \emptyset_{W_n}^{+}    \right\}.  $$
The rest of cascades we leaved unchanged. We denote obtained cascase $\widetilde{W}_n$. 
\vspace{0.5cm}

Put $U_0 = \omega$. Assume that $U_0, U_1, \ldots , U_{n-1}$ was defined, but it is impossible to define $U_n$. This means that every set $U\in \int \widetilde{W}_n $ is contained in  $\langle \bigcup_{i< n } {\cal V}_i \rangle$. 
On the other site $\max \widetilde{W}_n \in {\cal W}$ and so the family
$\{ U\cap  \max \widetilde{W}_n : U \in \bigcup_{i\leq n} {\cal V}_i \} $ has finite intersection property. By theorem of Dolecki $\{ U\cap  \max \widetilde{W}_n : U \in \bigcup_{i\leq n} {\cal V}_i \} $ is not contained in the contour  $\int \widetilde{W}_n$ of rank $\alpha_n +2$. A contradiction.
On each step of induction we can put $\bigcap_{i\leq n} U_i$ instead of $U_n$ and assume that the sequence $(U_n)_{n<\omega}$ is decreasing. 
\vspace{0.5cm}

Let $U= \bigcap_{n<\omega} U_n$. Conditions (1)-(4) guarantee that $U\in \int W$ and
$U \notin \langle \bigcup_{n<\omega} {\cal V}_n \rangle$. Indeed, assume that $U \in \langle \bigcup_{n<\omega} {\cal V}_n \rangle$, then there is a finite $M<\omega$ such that 
$U \in \langle \bigcup_{n<M} {\cal V}_n \rangle$. But $U_M  \notin \bigcup_{n\leq M} {\cal V}_n \rangle$ and $U \subset U_M$. Thus $U \notin \langle \bigcup_{n<\omega} {\cal V}_n \rangle$. A contradiction.

 \hspace{-2mm}
\gd

\medskip

{\small\sc \noindent {Micha\l} Machura, Instytut Matematyki, Uniwersytet {\'{S}l\c{a}ski},
, Katowice, Poland

E-mail:  machura@math.us.edu.pl}

{\small\sc \noindent Andrzej Starosolski, {Wydzia\l}  {Matematyki Stosowanej},
Politechnika \'{S}l\c{a}ska, Gliwice, Poland

E-mail:  andrzej.starosolski@polsl.pl}

\end{document}